\documentclass[12pt]{amsart}
\newtheorem{theorem}{Theorem}
\newtheorem{proposition}[theorem]{Proposition}
\newtheorem{corollary}[theorem]{Corollary}
\newtheorem{lemma}[theorem]{Lemma}
\def\ds{\displaystyle}

\title[The Kobayashi-Buseman pseudometric]{On the definition of the
Kobayashi-Buseman pseudometric }

\author{Nikolai Nikolov and Peter Pflug}

\address
{Institute of Mathematics and Informatics\\ Bulgarian Academy of
Sciences\\ Acad. G. Bonchev 8, 1113 Sofia,
Bulgaria}\email{nik@math.bas.bg}

\address{Carl von Ossietzky Universit\"at Oldenburg\\
Fachbereich Mathematik\\ Postfach 2503\\ D-26111 Oldenburg,
Germany}\email{pflug@mathematik.uni-oldenburg.de}



\begin{document}

\begin{thanks}{This note was written during the stay of the first
named author at the Universit\"at Oldenburg supported by a grant
from the DFG (January -- March 2006). He likes to thank both
institutions for their support.}
\end{thanks}

\begin{abstract} We prove that the $(2n-1)$-th Kobayashi pseudometric
of any domain $D\subset\Bbb C^n$ coincides with the
Kobayashi--Buseman pseudometric of $D,$ and that $2n-1$ is the
optimal number, in general.
\end{abstract}

\maketitle

\section{Introduction and results}

Let $\Bbb D\subset\Bbb C$ be the unit disc. Recall first the
definitions of the Lempert function $\tilde k_D$ and the
Kobayashi--Royden pseudometric $k_D$ of a domain $D\subset\Bbb
C^n$ (cf. \cite{Jar-Pfl}):
$$\aligned
\tilde
k_D(z,w)&=\inf\{\tanh^{-1}|\alpha|:\exists\varphi\in\mathcal
O(\Bbb D,D): \varphi(0)=z,\varphi(\alpha)=w\},\\
\kappa_D(z;X)&=\inf\{\alpha\ge 0:\exists\varphi\in\mathcal O(\Bbb
D,D): \varphi(0)=z,\alpha\varphi'(0)=X\},
\endaligned$$
where $z,w\in D,$ $X\in\Bbb C^n.$ The Kobayashi pseudodistance
$k_D$ can be defined as the largest pseudodistance which does not
exceed $\tilde k_D.$ Note that if $k_D^{(m)}$ denotes the $m$-th
Lempert function of $D$, that is,
$$k_D^{(m)}(z,w)=\inf\{\sum_{j=1}^m\tilde
k_D(z_j,z_{j+1}):z_1,\dots,z_m\in D, z_1=z,z_m=w\},$$ then
$$k_D(z,w)=\inf_m
k_D^{(m)}(z,w)=\inf\{\int_0^1\kappa_D(\gamma(t);\gamma'(t))dt\},$$
where the infimum is taken over all piecewise $C^1$-curves
$\gamma:[0,1]\to D$ connecting $z$ and $w$. By a result of M.~Y.
Pang (see \cite{Pang}), the Kobayashi--Royden pseudometric is the
infinitesimal form of the Lempert function for taut domains; more
precisely, if $D$ is a taut domain, then
$$\ds\kappa_D(z;X)=\lim_{\Bbb
C_\ast\ni t\to 0}\frac{\tilde k_D(z,z+tX)}{t}.\leqno{(1)}$$

In \cite{KobS}, S.~Kobayashi introduces a new invariant
pseudometric, called the Kobayashi--Buseman pseudometric in
\cite{Jar-Pfl}. One of the equivalent ways to define the
Kobayashi--Buseman pseudometric $\hat\kappa_D$ of $D$ is just to
set $\hat\kappa_D(z;\cdot)$ to be largest pseudonorm  which does
not exceed $\kappa_D(z;\cdot).$ Recall that
$$
\hat\kappa_D(z;X)=\inf\{\sum_{j=1}^m\kappa_D(z;X_j):m\in\Bbb N,\
\sum_{j=1}^mX_j=X\}.
$$ Thus, it is natural to consider the new function
$\kappa_D^{(m)},$ namely,
$$
\kappa_D^{(m)}(z;X)=\inf\{\sum_{j=1}^m\kappa_D(z;X_j):
\sum_{j=1}^mX_j=X\}.
$$
We call $\kappa_D^{(m)}$ the $m$-th Kobayashi pseudometric of $D.$
It is clear that $\kappa_D^{(m)}\ge\kappa_D^{(m+1)}$ and if
$\kappa_D^{(m)}(z;\cdot)=\kappa_D^{(m+1)}(z;\cdot)$ for some $m,$
then $\kappa_D^{(m)}(z;\cdot)=\kappa_D^{(j)}(z;\cdot)$ for any
$j>m.$ It is shown in \cite{KobS} that
$$
\kappa_D^{(2n)}=\hat\kappa_D.\leqno{(2)}
$$

Let now $D\subset\Bbb C^n$ be a taut domain. We point out that,
using the equalities (1) and (2), M.~Kobayashi (see \cite{KobM})
shows that
$$
\hat\kappa_D(z;X)=\lim_{\Bbb C_\ast\ni t\to 0}\frac{
k_D(z,z+tX)}{t}.
$$
Obvious modifications in the proof of this result lead to
$$
\lim_{u,v\to z,\ u\neq
v}\frac{k_D^{(m)}(u,v)-\kappa_D^{(m)}(z;u-v)}{||u-v||}=0.\leqno{(3)}
$$
uniformly in $m$ and locally uniformly in $z;$ thus,
$$
\kappa_D^{(m)}(z;X)=\lim_{\Bbb C_\ast\ni t\to 0}\frac{
k_D^{(m)}(z,z+tX)}{t}
$$
uniformly in $m$ and locally uniformly in $z$ and $X.$

The aim of this note is the following result which improves (2).

\begin{theorem} For any domain $D\subset\Bbb C^n$ we have that
$$\kappa_D^{(2n-1)}=\hat\kappa_D.\leqno{(4)}$$

On the other hand, if $n\ge 2$ and $$D_n=\{z\in\Bbb
C^n:\sum_{j=2}^n(2|z_1^3-z_j^3|+|z_1^3+z_j^3|)<2(n-1)\},$$ then
$$\kappa_{D_n}^{(2n-2)}(0;\cdot)\neq\hat\kappa_{D_n}(0;\cdot).\leqno{(5)}$$
\end{theorem}

Note that the proof below shows that the equality (4) remains true
for any $n$-dimensional complex manifold.

An immediately consequence of Theorem 1 and the equality (3) is:

\begin{corollary} For any taut domain $D\subset\Bbb C^n$ one has
that
$$
\lim_{w\to z,w\neq z}\frac{k_D^{(2n-1)}(z,w)}{k_D(z,w)}=1
$$
locally uniformly in $z,$ and $2n-1$ is the optimal number, in
general.
\end{corollary}

\noindent{\bf Remarks.} (i) If $D\subset\Bbb C,$ then even $\tilde
k_D=k_D$ (cf. \cite{Jar-Pfl}).

(ii) Corollary 2 holds for $n$-dimensional taut complex manifolds.

(iii) Observe that Corollary 2 may be taken as a very weak version
of the following question asked by S.~Krantz (see \cite{Kra}):
whether there is a positive integer $m=m(D)$ such that
$k_D=k_D^{(m)}.$
\medskip

Let now $h_S$ be the Minkowski functions of a starlike domain
$S\subset\Bbb R^N,$ that is, $h_S(X)=\inf\{t>0:X/t\in S\}$). We
may define as above
$$
h_S^{(m)}(X)=\inf\{\sum_{j=1}^mh_S(X_j): \sum_{j=1}^mX_j=X\}.
$$
Then the Minkowski function $h_{\hat S}$ of the convex hull $\hat
S$ of $S$ is the largest pseudonorm  which does not exceed $h_S.$
It follows by a lemma due to C.~Carath\'eodory (cf. \cite{KobM})
that
$$
\aligned h_{\hat S}=h_S^{(N)}=\inf\{\sum_{j=1}^Mh_S(X_j):M\le N,\
\sum_{j=1}^MX_j=X,\\X_1,\dots, X_M\mbox{ are }\Bbb
R\mbox{-linearly independent}\}.\endaligned\leqno{(6)}
$$
One can easily see that $N$ is the optimal number for the class of
starlike domains in $\Bbb R^N.$

Denote by $I_{D,z}$ the indicatrix of $\kappa_D(z;\cdot)$, that
is, $I_{D,z}=\{X\in\Bbb C^n:\kappa_D(z;X)<1\}.$ Note that
$I_{D,z}$ is a balanced domain (a domain $B\subset\Bbb C^n$ is
said to be balanced if $\lambda X\in B$ for any
$\lambda\in\overline{\Bbb D}$ and any $X\in B$). In particular,
$I_{D,z}$ is a starlike domain and hence (2) follows by (6).
Similarly, (4) will follow by the following.

\begin{proposition} If $B\subset\Bbb C^n$ is a balanced domain,
then
$$h_{\hat B}=h_B^{(2n-1)}.\leqno{(7)}$$
\end{proposition}

Observe that the domain $D_n$ from Theorem 1 is pseudoconvex and
balanced, thus $\kappa_{D_n}(0;\cdot)=h_{D_n}$ (cf.
\cite{Jar-Pfl}) and so
$\kappa_{D_n}^{(m)}(0;\cdot)=h_{D_n}^{(m)}.$ Then inequality (5)
is equivalent to
$$
h_{\hat D_n}\neq h_{D_n}^{(2n-2)}.\leqno{(8)}
$$

\section{Proofs}

To prove Proposition 3, we shall need the following result.

\begin{lemma} Any balanced domain can be exhausted by bounded balanced
domains with continuous Minkowski functions.
\end{lemma}

\begin{proof} Let $B\subset\Bbb C^n$ be a balanced domain.
Denote by $\Bbb B_n(z,r)\subset\Bbb C^n$ the ball with center $z$
and radius $r.$ For $z\in\Bbb C^n$ and $j\in\Bbb N$, set
$F_{n,j,z}:=\overline{\Bbb B_n(z,||z||^2/j)}.$ We may assume that
$\Bbb B_n(0,1)\subset\subset B.$ Put
$$
B_j:=\{z\in \Bbb B_n(0,j): F_{n,j,z}\subset B\},\quad j\in\Bbb N.
$$
Then $(B_j)_{j\in\Bbb N}$ is an exhaustion of $B$ by non-empty
bounded open sets. We shall show that $B_j$ is a balanced domain
with continuous Minkowski functions.

For this, take any $z\in B_j$ and $0\neq\lambda\in\overline{\Bbb
D},$ and observe that $F_{n,j,\lambda z} \subset\lambda
F_{n,j,z}\subset B.$ Thus, $B_j$ is a balanced domain.

Since $h_{B_j}$ is an upper semicontinuous function, it remains to
prove that it is lower semicontinuous. Assuming the contrary, we
may find a sequence of points $(z_k)_{k\in\Bbb N}$ converging to
some point $z\in\Bbb C^n$ and a positive number $c$ such that
$h_{B_j}(z_k)<1/c<h_{B_j}(z)$ for any $k.$  Note that
$F_{n,j,cz_k}\subset B$, $k\in\Bbb N$. Hence $\Bbb
B_n(cz,c^2\|z\|^2/j)\subset B$. On the other hand, fix $t\in
(0,1)$ such that $h_{B_j}(tcz)>1$. Then $F_{n,j,tcz} \subset\Bbb
B_n(cz,c^2\|z\|^2/j)\subset B$; thus $h_{B_j}(tcz)<1$, a
contradiction.
\end{proof}

\noindent{\it Proof of Proposition 3.} First, we shall prove (7)
in the case, when $B\subset\Bbb C^n$ is a bounded balanced domain
with continuous Minkowski function. Fix a vector $X\in\Bbb
C^n\setminus\{0\}.$ Then  $h_{\hat B}(X)\neq 0$ and we may assume
that $h_{\hat B}(X)=1.$ By the continuity of $h_B$ and (6), there
exist $\Bbb R$-linearly independent vectors $X_1,\dots,X_m$ ($m\le
2n$) such that $\ds\sum_{j=1}^m X_j=X$ and
$\ds\sum_{j=1}^mh_B(X_j)=1.$ Since $h_{\hat B}$ is a norm, the
triangle inequality implies that $h_B(X_j)=h_{\hat B}(X_j),$\;
$j=1,\dots,m.$ To prove (7), it suffices to show that $m\neq 2n.$
The convexity of $\hat B$ provides a support hyperplane $H$ for
$\hat B$ at $X\in\partial\hat B,$ say $H=\{z\in\Bbb
C^n:\mbox{Re}\langle z-X,\overline X_0\rangle=0\},$ $X_0\in\Bbb
C^n,$ where $\langle\cdot,\cdot\rangle$ stands for the Hermitian
scalar product in $\Bbb C^n.$ Assuming $m=2n$ implies that $\ds
H=\{\sum_{j=1}^m\alpha_j \hat X_j:\sum_{j=1}^m\alpha_j =1,\
\alpha_1,\dots,\alpha_m\in\Bbb R\},$ where $\hat
X_j:=X_j/h_B(X_j)\in\partial\hat B$. In particular, $\partial\hat
B$ contains a relatively open subset of $H.$ Since $\hat B$ is a
balanced domain, it follows that its intersection with the plane,
spanned by $X_0,$ is a disc whose boundary contains a line
segment, a contradiction.

Now let $B\subset\Bbb C^n$ be an arbitrary balanced domain. If
$(B_j)_{j=1}^\infty$ is an exhaustion of $B$ given by Lemma 4,
then $h_{B_j}\searrow h_B$ pointwise and hence $h_{\hat
B_j}\searrow h_{\hat B}$ by (6). Then (7) follows by the
inequalities $h_{\hat B}\le h_B^{(2n-1)}\le h_{B_j}^{(2n-1)}$ and
the equality $h_{\hat B_j}=h_{B_j}^{(2n-1)}$ from above.\qed
\medskip

\noindent{\it Proof of the inequality (8).} Let $L_n=\{z\in \Bbb
C^n:z_1=1\}.$ Then the triangle inequality implies that
$D_n\subset\Bbb D\times\Bbb C^{n-1}$ and
$$
F_n:=\partial D_n\cap L_n=\{z\in\Bbb C^n:z_1=1,\ z_j\in\Omega,\
2\le j\le n\},
$$
where $\Omega$ is the set of the third roots of unity. Denoting by
$\Delta$ the convex hull of $\Omega,$ it follows that
$$
\partial\hat D_n\cap L_n=\hat F_n=\{1\}\times\Delta^{n-1}.
$$
Hence, $\partial\hat D_n\cap L_n$ is a $(2n-2)$-dimensional convex
set. Put $\tilde F_n=\{Y\in\hat F_n:h_{D_n}^{(2n-2)}(Y)=1\}.$ If
$X\in\tilde F_n,$ then there exist $X_1,\dots,X_m\in\Bbb
C^n\setminus\{0\},$ $m\le 2n-2,$ such that $\ds\sum_{j=1}^m X_j=X$
and $\ds\sum_{j=1}^m h_{D_n}(X_j)=1$ (note that $D_n$ is taut).
Hence, $X_1/h_{D_n}(X_1),\dots,X_m/h_{D_n}(X_m)\in F_n$ and $X$
belongs to their convex hull. Since $F_n$ is a finite set, it
follows that $\tilde F_n$ is contained in a finite union of at
most $(2n-3)$-dimensional convex sets. Thus, $\hat F_n\neq\tilde
F_n$ which implies that $h_{\hat D_n}\neq h_{D_n}^{(2n-2)}.$\qed

\end{document}